\documentclass{article}%
\usepackage{amsmath}%
\usepackage{amsfonts}%
\usepackage{amssymb}%
\usepackage{graphicx}

\begin{document}

\title{A Closed-Form Solution to the Arbitrary Order Cauchy Problem with Propagators}
\author{Henrik Stenlund\thanks{The author is grateful to Visilab Signal Technologies for supporting this work.}
\\Visilab Signal Technologies Oy, Finland}
\date{November 24, 2014}
\maketitle

\begin{abstract}
The general abstract arbitrary order (N) Cauchy problem was solved in a closed form as a sum of exponential propagator functions. The infinite sparse exponential series was solved with the aid of a homogeneous differential equation. It generated a linear combination of exponential functions. The Cauchy problem solution was formed with N linear combinations of N exponential propagators. \footnote{Visilab Report \#2014-12} \footnote{Dedicated to Peter, my younger brother}
\subsection{Keywords}
arbitrary order Cauchy problem, differential equations, initial value problems, infinite series, summation of series
\subsection{Mathematical Classification}
Mathematics Subject Classification 2010: 47D09, 34A12, 35E15, 35F10, 35F40, 35G10, 35G40, 44A10, 11M41, 16W60, 20F14, 40A25, 65B10
\end{abstract}

\section{Introduction}
The abstract general arbitrary order $N$ Cauchy problem ($CP_N$) is the following differential equation, treated as an initial value problem. 
\begin{equation}
\frac{\partial^N{u(x,t)}}{\partial{t^N}}=G{u(x,t)}  \label{eqn500}
\end{equation}
The initial value functions are defined as follows ($i=0,1,2...N-1$)
\begin{equation}
[\frac{\partial^i{u(x,t)}}{\partial{t^i}}]_{t=t_0}=u_i(x)  \label{eqn510}
\end{equation}
$G$ is a linear operator independent of $t$ and $\frac{\partial}{\partial{t}}$
\begin{equation}
G=G(x,\frac{\partial}{\partial{x}})  \label{eqn520}
\end{equation}
This is a general abstract problem with ($t \in{R},x \in{R}$). We assume $u\in{R}$ but it may adopt complex values. It was recently solved by the author \cite{Stenlund2014} with an infinite propagator series. Well-posedness is assumed for the ($CP_N$) and that subject has been treated thoroughly, for example in \cite{Fattorini1983}, \cite{Melnikova2001} and \cite{Xiao1998} having excellent references to earlier work. We assume that $u(x,t)$ can be expanded as a Taylor's power series,  requiring:

\itshape
- $u(x,t)$ is continuous as a function of $t, t > t_0$ 

- the series converges

- all the derivatives exist for any $n=1,2,3..$:
\begin{equation}
\left|\frac{\partial{^{n}u(x,t)}}{\partial{t^n}}\right| < \infty  \nonumber
\end{equation}
\normalfont
We could equally well treat it as an analytic function on the complex plane with the same results. The Taylor's series is
\begin{equation}
u(x,t)=\sum^{\infty}_{n=0}{\frac{(t-t_0)^{n}}{n!}[\frac{\partial{^{n}u(x,t)}}{\partial{t^n}}]_{t=t_0}} \label{eqn525}
\end{equation}
and we can solve for all the derivatives in a cyclic manner
\begin{equation}
\frac{\partial{^{n}u(x,t)}}{\partial{t^n}}= \nonumber
\end{equation}
\begin{equation}
{G^{\frac{n}{N}}u(x,t), \ \ \ n=0,N,2N,3N...}  \nonumber
\end{equation}
\begin{equation}
{G^{\frac{n-1}{N}}\frac{\partial{u(x,t)}}{\partial{t}}, \ \ \ n=1,N+1,2N+1,3N+1...}  \nonumber
\end{equation}
\begin{equation}
{G^{\frac{n-2}{N}}\frac{\partial{^{2}u(x,t)}}{\partial{t^2}}, \ \ \ n=2,N+2,2N+2,3N+2...}  \nonumber
\end{equation}
...
\begin{equation}
{G^{\frac{n-N+1}{N}}\frac{\partial{^{N-1}u(x,t)}}{\partial{t^{N-1}}}, \ \ \ n=N-1,N+N-1,2N+N-1...} \label{eqn540}
\end{equation}
The ($CP_N$) solution will be
\begin{equation}
u(x,t)=\sum^{N-1}_{j=0}{(t-t_0)^j\sum^{\infty}_{n=0}{\frac{((t-t_0)^{N}G)^n}{(Nn+j)!}}}u_j(x), \ \ \ t > t_0   \label{eqn580}
\end{equation}
The latter summation is an infinite series of sparse exponential type for which there is seems not to be available any solution in a closed form. The motivation for this paper is to find a simplified solution to the ($CP_N$). 

In Chapter 2 we solve the sparse exponential power series. In Chapter 3 we show the new solution to the ($CP_N$). In Chapter 4 we present a simple example application of the resulting solution. We try to keep the presentation simple.

\section{Solving the Infinite Sparse Exponential Series}
We start by studying the sparse exponential series which is defined as follows
\begin{equation}
y_j(z)=\sum^{\infty}_{k=j,N+j,2N+j..}{\frac{z^k}{k!}}   \label{eqn2000}
\end{equation}
where $N\in{R^+}$, $z\in{R}$ and $j=0..N-1$. The $j$ and $N$ are external parameters depending on the particular problem setup. $z$ can be complex as well. The series resembles an exponential function's power series very closely but is missing terms. The series converges absolutely just like the exponential function's series does. 

Regular summation methods do not offer much help for this series having a factorial denominator. We have to improvise a new method as follows. We differentiate (\ref{eqn2000}) $N$ times to get
\begin{equation}
\frac{\partial^{N}y_j(z)}{\partial{z^N}}=\sum^{\infty}_{N+j,2N+j,3N+j..}{\frac{z^{k-N}}{(k-N)!}}=\sum^{\infty}_{m=j,N+j,2N+j..}{\frac{z^m}{m!}}  \label{eqn2030}
\end{equation}
where we changed the index to $k-N=m$. We have thus returned to the original series. This is apparently some sort of an exponential function behaving as follows 
\begin{equation}
\frac{\partial^{N}y_j(z)}{\partial{z^N}}=y_j(z)  \label{eqn2040}
\end{equation}
It is an $N$th order homogeneous differential equation and the characteristic equation is
\begin{equation}
D^{N}-1=0  \label{eqn2065}
\end{equation}
The roots are
\begin{equation}
(\lambda_{n,j})^{N}=1=e^{{2\pi{in}}}  \label{eqn2060}
\end{equation}
and
\begin{equation}
\lambda_{n,j}=e^{\frac{2\pi{in}}{N}}  \label{eqn2070}
\end{equation}
where $n\in{R^+}$. We can drop the $j$ since $\lambda$ does not carry any dependence of it. The solution to the differential equation (\ref{eqn2040}) is 
\begin{equation}
y_j(z)=\sum^{N}_{n=1}{C_{n,j}{e^{\lambda_{n}{z}}}}\equiv{\sum^{\infty}_{k=j,N+j,2N+j..}{\frac{z^k}{k!}}}  \label{eqn2080}
\end{equation}
The solution is a sum of linearly independent exponential functions. We are left with the constant coefficients $C_{n,j}$ as unknowns. They are solved by applying obvious initial conditions. We evaluate the identity (\ref{eqn2080}) and its derivatives at $z=0$. Thus we get the following system of equations with the aid of Kronecker delta functions for a fixed $j$
\begin{equation}
y_j(0)=\sum^{N}_{n=1}{C_{n,j}\cdot{({\lambda_{n}})^0}}=\delta_{j,0}   \label{eqn2090}
\end{equation}
\begin{equation}
y_j^{(1)}(0)=\sum^{N}_{n=1}{C_{n,j}\cdot{({\lambda_{n}})^1}}=\delta_{j,1}   \label{eqn2100}
\end{equation}
\begin{equation}
y_j^{(2)}(0)=\sum^{N}_{n=1}{C_{n,j}\cdot{({\lambda_{n}})^2}}=\delta_{j,2}   \label{eqn2110}
\end{equation}
\begin{equation}
...
\end{equation}
\begin{equation}
y_j^{(N-1)}(0)=\sum^{N}_{n=1}{C_{n,j}\cdot{({\lambda_{n}})^{N-1}}}=\delta_{j,{N-1}}   \label{eqn2120}
\end{equation}
The right side equalities form a system of linear equations for the $C_{n,j}$ for a particular value of $j$. For our ($CP_N$) we have $N$ systems of equations or $N^2$ equations and the number of coefficients $C_{n,j}$ is $N^2$. The systems are separate since the $C_{n,j}$'s have a different $j$ in each system and can thus be solved independently in a straightforward way with matrices as usual, assuming there are solutions. We are very close to the Wronskian determinant of the set of exponential functions, evaluated at $z=0$. In lowest order cases we have for $N=1$
\begin{equation}
C_{1,0}=1
\end{equation}
which will lead to the regular exponential function solution. For $N=2$ we have
\begin{equation}
C_{1,0}=\frac{1}{2}
\end{equation}
\begin{equation}
C_{2,0}=\frac{1}{2}
\end{equation}
\begin{equation}
C_{1,1}=\frac{-1}{2}
\end{equation}
\begin{equation}
C_{2,1}=\frac{1}{2}
\end{equation}
producing $cosh()$ and $sinh()$ functions. Higher order solutions become much more complicated but can be facilitated either by symbolic or numeric processing. The solution (\ref{eqn2080}) is independent of the ($CP_N$). 
\section{Solution of the Arbitrary Order Cauchy Problem in a Closed Form}
We need to process our equation (\ref{eqn580}) a little so that we can apply the results in the preceding chapter. If we substitute an operator
\begin{equation}
\kappa=(t-t_0)G^{\frac{1}{N}}
\end{equation}
the equation transforms to
\begin{equation}
u(x,t)=\sum^{N-1}_{j=0}{\frac{(t-t_0)^j}{\kappa^{j}}\sum^{\infty}_{n=0}{\frac{\kappa^{Nn+j}}{(Nn+j)!}}}u_j(x)  \label{eqn3000}
\end{equation}
Next we change the index to $k=Nn+j$ and we get
\begin{equation}
u(x,t)=\sum^{N-1}_{j=0}{\frac{1}{G^{\frac{j}{N}}}\sum^{\infty}_{k=j,N+j,2N+j,3N+j..}{\frac{\kappa^{k}}{k!}}}u_j(x)   \label{eqn3030}
\end{equation}
We have the latter sum ready for applying the previous result (\ref{eqn2080}) with $z=\kappa$ and substitute to obtain
\begin{equation}
u(x,t)=\sum^{N-1}_{j=0}{\frac{1}{G^{\frac{j}{N}}}\sum^{N}_{n=1}{C_{n,j}e^{\lambda_n{(t-t_0)G^{\frac{1}{N}}}}}}u_j(x), \ \ \ t > t_0   \label{eqn3040}
\end{equation}
where
\begin{equation}
\lambda_{n}=e^{\frac{2\pi{in}}{N}}  \label{eqn3070}
\end{equation}
The $C_{n,j}$ are the coefficients to be solved and they are determined by the $N$ alone, having no other connection to the ($CP_N$). We have obtained a general closed-form solution to the ($CP_N$). It is interesting that it always consists of a combination of linearly independent exponential functions or propagators. On the other hand, the $G$ does not dictate the form of the solution at this level.

\section{An Example}
\subsection{The Nth Order Wave Equation}
We apply the results above to the problem
\begin{equation}
G=v^N{\frac{\partial^N}{\partial{x^N}}}  \label{eqn12030}
\end{equation}
The initial value functions are defined as follows ($i=0,1...N-1$)
\begin{equation}
[\frac{\partial^i{u(x,t)}}{\partial{t^i}}]_{t=t_0}=u_i(x)  \label{eqn12040}
\end{equation}
and the ($CP_N$) is
\begin{equation}
\frac{\partial^N{u(x,t)}}{\partial{t^N}}=v^N{\frac{\partial^N}{\partial{x^N}}}{u(x,t)}  \label{eqn12050}
\end{equation}
We use the solution (\ref{eqn3040}) and substitute, getting
\begin{equation}
u(x,t)=\sum^{N-1}_{j=0}{\frac{1}{(\frac{{\delta}^{\frac{j}{N}}}{\delta{x}^{\frac{j}{N}}})}\sum^{N}_{n=1}{C_{n,j}e^{v(t-t_0)\lambda_n{\frac{\delta}{\delta{x}}}}}}u_j(x)   \label{eqn12055}
\end{equation}
Since the exponential operator is a universal propagator,
\begin{equation}
u(x,t)=\sum^{N-1}_{j=0}{\frac{1}{D^{\frac{j}{N}}}\sum^{N}_{n=1}{C_{n,j}}}u_j(x+v(t-t_0)\lambda_n)   \label{eqn12058}
\end{equation}
The operator in the first summation is an integration operator of $j/N$ times and may be fractional. This solution consists of $N^2$ waves traveling to $N$ directions from the origin on the complex plane. 
\section{Discussion}
We have presented a solution to the abstract general arbitrary order Cauchy problem ($CP_N$). It was recently solved as an infinite series of propagators, equation (\ref{eqn580}). 

The infinite summation part, which is a sparse exponential power series, is identified as a function having an $N$th order homogeneous differential equation. That can be solved with a linear combination of $N$ exponential functions which are linearly independent of each other, equation (\ref{eqn2080}). The coefficients are solved from a simple linear system of equations. Since the infinite ($CP_N$) solution consists of two summations, the outer summation having N terms. Therefore, the ($CP_N$) solution will consist of $N$ sets of $N$ linear combinations of exponential propagators or $N^2$ of propagators. 

The solution (\ref{eqn3040}) will, at least superficially, contain fractional operators. The infinite solution (\ref{eqn580}) would indicate that there are no fractional operators afterall. 

The main results of this work are equations (\ref{eqn2080}) and (\ref{eqn3040}). As a byproduct in this study we have found a general solution to the sparse exponential series. 


\end{document}